\documentclass[conference]{IEEEtran}
\IEEEoverridecommandlockouts % allows use of \thanks command

\usepackage{caption}
\usepackage{graphicx}
\usepackage{amssymb,amsmath,nicefrac,amsthm,multicol, mathtools}
\usepackage{color}
\usepackage{tocvsec2}
\usepackage{cite}
\usepackage{booktabs} % for nice tables
\usepackage{bm}
\usepackage{multirow}
\usepackage{epstopdf} % solve the problem with Figures
\usepackage{setspace}
\usepackage{enumerate}
\usepackage[utf8]{inputenc}
\usepackage[T1]{fontenc}
\usepackage{url}
\usepackage{subfigure}
\usepackage{mathtools}
\usepackage{derivative}
\usepackage[scr=boondox]  % heavily sloped
           {mathalpha}
\usepackage{tabularx}

\newcommand{\oft}{}

\newcommand{\J}{J(\tilde{x}_{k:k+N}, \tilde{u}_{k:k+N-1})}
\newcommand{\fd}{f_d(\tilde{x}_{k+i}, \tilde{u}_{k+i}, \tilde{\theta})}
\newcommand{\F}{\mathcal{F}(x\oft)}
\newcommand{\G}{\mathcal{G}(x\oft, u\oft)}

\newcommand{\overm}[1]{\frac{#1}{\mu}}

\newcommand{\threeDvec}[3]{
    \begin{bmatrix}
        #1 \\ #2 \\ #3
    \end{bmatrix}
}
\newcommand{\threeDvecT}[3]{
    \begin{bmatrix}
        #1 & #2 & #3
    \end{bmatrix}^T
}
\newcommand{\fourDvec}[4]{
    \begin{bmatrix}
        #1 \\ #2 \\ #3 \\ #4
    \end{bmatrix}
}
\newcommand{\fourDvecT}[4]{
    \begin{bmatrix}
        #1 & #2 & #3 & #4
    \end{bmatrix}^T
}

\usepackage{mathtools}
\begin{document}
\title{Efficient Estimation of Relaxed Model Parameters for Robust UAV Trajectory Optimization}

\author{\IEEEauthorblockN{Derek Fan}
\IEEEauthorblockA{\textit{Department of Mechanical Engineering} \\
\textit{Carnegie Mellon University}\\
Pittsburgh, PA, USA \\
derekfan@andrew.cmu.edu}
\and
\IEEEauthorblockN{David A. Copp}
\IEEEauthorblockA{\textit{Department of Mechanical and Aerospace Engineering} \\
\textit{University of California, Irvine}\\
Irvine, CA, USA \\
dcopp@uci.edu}
}
\maketitle

\begin{abstract}
    Online trajectory optimization and optimal control methods are crucial for enabling sustainable unmanned aerial vehicle (UAV) services, such as agriculture, environmental monitoring, and transportation, where available actuation and energy are limited. However, optimal controllers are highly sensitive to model mismatch, which can occur due to loaded equipment, packages to be delivered, or pre-existing variability in fundamental structural and thrust-related parameters. To circumvent this problem, optimal controllers can be paired with parameter estimators to improve their trajectory planning performance and perform adaptive control. However, UAV platforms are limited in terms of onboard processing power, oftentimes making nonlinear parameter estimation too computationally expensive to consider. To address these issues, we propose a relaxed, affine-in-parameters multirotor model along with an efficient optimal parameter estimator. We convexify the nominal Moving Horizon Parameter Estimation (MHPE) problem into a linear-quadratic form (LQ-MHPE) via an affine-in-parameter relaxation on the nonlinear dynamics, resulting in fast quadratic programs (QPs) that facilitate adaptive Model Predictve Control (MPC). This makes real-time applications more feasible. We compare this approach to the equivalent nonlinear estimator in Monte Carlo simulations, demonstrating a decrease in average solve time by \textbf{98.2\%} and trajectory optimality cost by \textbf{23.9-56.2\%}.
\end{abstract}


\section{Introduction}\label{sec:intro}

Unmanned aerial vehicles (UAVs) are rapidly transforming our society. Their autonomous capabilities have opened up possibilities for a wide range of applications \cite{https://doi.org/10.1002/rob.21962}, including many focused on improving sustainability in transportation, agriculture, and environmental monitoring. Some of these applications include efficient package delivery \cite{chiang2019impact, saunders2024autonomous}, wildfire monitoring \cite{bailon2022real}, and aerial crop monitoring and spraying \cite{radoglou2020compilation}. These UAV applications benefit from trajectory optimization and optimal control methods such as Model Predictve Control (MPC) \cite{mayne2000constrained} for efficient autonomous behaviors. These algorithms are powerful tools for controlling challenging, underactuated systems (e.g., UAVs), optimizing between objectives such as actuator effort and energy consumption \cite{doi:10.2514/1.G002507}, and continuously re-planning around disturbances and obstacles. The value of MPC is especially clear when considering an application such as UAV package delivery, which has been found to be more energy-efficient than current truck delivery services \cite{RODRIGUES2022100569, KOIWANIT2018201}.

However, for MPC to demonstrate optimal behavior, the dynamics model that it employs must accurately reflect the system that is being controlled. Since MPC is sensitive to modeling errors, and performing system identification beforehand can be impractical, it may be necessary to implement modifications that improve its robustness to model uncertainty. Model mismatch can occur in drone delivery, where the system's aerodynamic and inertial properties dramatically change before and after shipment. Other applications include agricultural crop spraying and environmental monitoring, where various attachments (e.g., sprayers, cameras, other sensors) introduce uncertainty or variability in model parameters. Although this issue can be addressed by pairing MPC with a parameter estimator (see, e.g., \cite{copp2016addressing}), nonlinearities and limited onboard computing power make online estimation challenging.

Modeling nonlinear systems with linear dynamics enables the usage of mature techniques from linear control, estimation, and optimization literature. The conventional strategy is to repeatedly linearize the nominal nonlinear dynamics about the system's evolving state. Although these linear approximations are oftentimes serviceable, the linearization process introduces additional algorithmic complexity and approximation error. Fortunately, if the system is nonlinear in parameters solely due to multiplicative coupling, a simple change of variables can relax the model without the need for iterative re-linearization \cite{9290199, doi:10.2514/1.G005376}.

For optimal parameter estimation with conventional quadratic penalties, linear models allow access to simple and effective methods such as the Kalman filter \cite{kalman1960new}. However, the Kalman filter is only capable of solving for unconstrained optimal estimates in the presence of Gaussian disturbances \cite{haseltine2005critical}. This is problematic when considering that environmental disturbances may be non-Gaussian, and that explicit constraints on a relaxed set of parameters are desirable for tight search spaces and dynamically feasible values. Moving Horizon Estimation (MHE) overcomes these weaknesses by employing constrained optimization, and it can handle non-Gaussian noise and disturbance distributions \cite{rao2001constrained,rao2003constrained}. It is often implemented as a nonlinear program (NLP) to address a general set of estimation problems. For example, the coordination of multiple UAVs for target tracking using MPC and MHE with nonlinear UAV dynamics was presented in \cite{quintero2015robust} and included estimation of unmeasured disturbances in \cite{quintero2017robust}.

In contrast to generalized NLPs, which can vary dramatically in terms of complexity and difficulty to solve, QPs are among the most well-understood class of optimization problems \cite{PANG1983583}. Their structure guarantees global optimality and is amenable to efficient solvers with strong convergence properties. Furthermore, there has been recent work in accelerating QP solve times on microcontrollers \cite{nguyen2024tinympc} and graphical processing units (GPUs) \cite{bishop2024relu}. These algorithms illustrate some of the numerical benefits of reformulating the parameter estimation problem as a QP.

We propose an optimal parameter estimator that facilitates fast adaptation and optimal control even under significantly uncertain UAV model parameters. Specifically, we reformulate the nonlinear multirotor UAV dynamics into a relaxed, affine-in-parameter model. Said model is implemented within a Moving Horizon Parameter Estimator (MHPE) to form a linear-quadratic variant (LQ-MHPE) that exploits the QP problem structure, enabling efficient parameter estimation for adaptive nonlinear model predictive control (NMPC).

To demonstrate the effectiveness of our proposed estimator, we run Monte Carlo simulations across two different UAV models with significantly varying model parameters, multiple numerical optimization solvers, and various randomized initial conditions. Specifically, we compare nonlinear MHPE (NMHPE) with the convexified LQ-MHPE and show that LQ-MHPE runs with faster solve times and adaptively updates the model used in the NMPC to produce comparable, and in some cases better, trajectory optimality. We thereby show that LQ-MHPE enables adaptive NMPC with sufficient speed and accuracy to enhance UAV performance in sustainability applications.

This paper is organized as follows. First, we formulate the generic system dynamics, trajectory optimization, and optimal estimation problem in Sections \ref{sec:prob-form}, \ref{sec:nmpc} and \ref{sec:parameter-estimation}, respectively. Next, we present the nominal multirotor dynamics and affine-in-parameter model derivation in Section \ref{sec:nl-dynamics}. We showcase the implementation details, simulation trial setup, and results in Section \ref{sec:results}, and conclude in Section \ref{sec:conclusion}.

\section{Problem Formulation}\label{sec:prob-form}
We seek to control and estimate the parameters of the following continuous-time nonlinear system that can represent UAV dynamics with potentially unknown parameter values:
\begin{align}\label{eq:generic-dynamics}
    \dot{x} = f(x\oft, u\oft, \theta) + w\oft,
\end{align}
where $x\oft \in \mathbb{R}^n $ is the state,
$u\oft \in \mathbb{R}^m$ is the control input,
and $w\oft \in \mathcal{W} \subset \mathbb{R}^n$ is additive Gaussian process noise at time $t \in \mathbb{R}$. In addition, $\theta \in \mathbb{R}^p$ is a vector of constant unknown parameters. We can discretize the dynamics in \eqref{eq:generic-dynamics} as
\begin{align}
    x_{k+1} = f_d(x_k, u_k, \theta) + w_k,
\end{align}
where $k \in \mathbb{N}$ is the discrete time step and $f_d$ is an arbitrary numerical integration of $f$.

\begin{figure}[t!]
    \includegraphics[width=\linewidth]{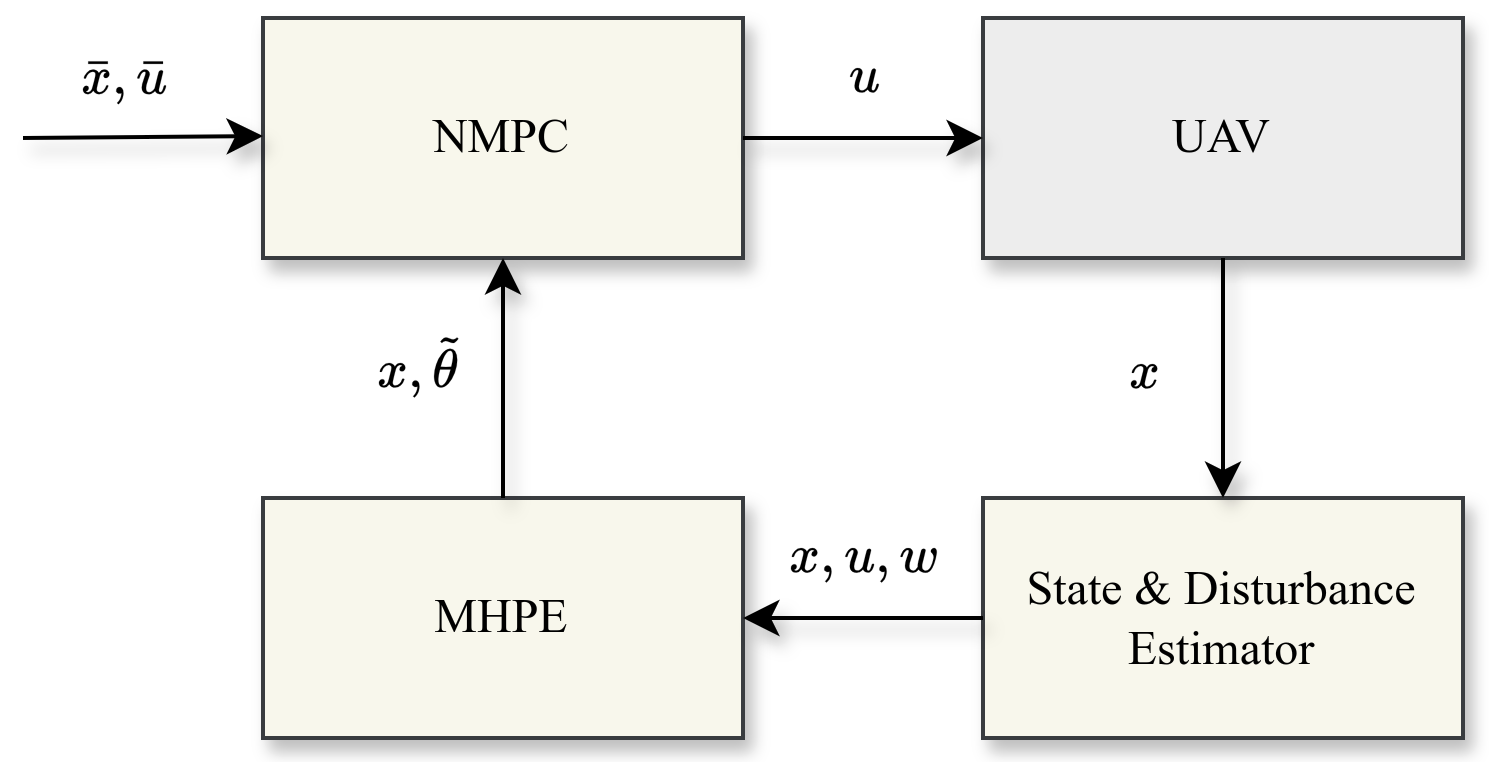}
    \caption{The flow of information within a control loop of MHPE-based adaptive NMPC. Although we assume information of the states and disturbances, a separate estimation filter would provide such estimates to the MHPE in practice.}
\end{figure}

\section{Nonlinear Model Predictve Control}\label{sec:nmpc}
NMPC is a powerful optimal control strategy for systems such as \eqref{eq:generic-dynamics} that optimizes performance metrics while explicitly handling constraints. For drone delivery, NMPC is especially effective for finding local energy-efficient trajectories.
At every time step, NMPC solves a trajectory optimization problem over a finite time horizon to find a sequence of future control inputs and applies the first input in said sequence. This trajectory optimization can be formulated as the following NLP:
\begin{subequations}\label{eq:nmpc}
\begin{align}
    \min_{\tilde{{x}}_{k:k+N}, \tilde{{u}}_{k:k+N-1}}
        && &\J,& \label{eq:J-generic} \\
    \text{s.t.}
        && &\tilde{{x}}_{k} = {x}_k,&\label{eq:initial-state} \\
        && &\tilde{{x}}_{k+i+1} = \fd,&\label{eq:eq-constr} \\
        && &\tilde{{x}}_{k:k+N} \in {\mathcal{X}},& \label{x-constr}\\
        && &\tilde{{u}}_{k:k+N-1} \in {\mathcal{U}},& \label{u-constr}
\end{align}
\end{subequations}
where $N \in \mathbb{N}$ is the length of the trajectory (prediction horizon), $i \in \{0,1,\hdots,N-1\}$ corresponds to a time step in the horizon, $\tilde{x}_{k:k+N}$ denotes the sequence of predicted states, $\tilde{u}_{k:k+N-1}$ denotes the sequence of predicted control inputs, and ${\tilde\theta}$ is an estimate of ${\theta}$ that is constant throughout each computed trajectory. Constraint \eqref{eq:initial-state} ensures that the initial predicted state matches the measured current state, \eqref{eq:eq-constr} ensures the predicted states follow the dynamics, and \eqref{x-constr} and \eqref{u-constr} define state and actuator constraints, respectively.

For an NMPC formulation that optimizes state trajectory tracking and energy consumption, one can choose the objective function $J(\cdot)$ in problem \eqref{eq:nmpc} as a quadratic function:
\begin{align}
\begin{split}\label{eq:J-quad}
    J(\tilde{x}_{k:k+N},& \tilde{u}_{k:k+N-1}) :=
        || \tilde{x}_{k+N} - \bar{x}_{k+N} ||^2_{{Q_f}} \\
    +& \sum^{N-1}_{i=0}
        || \tilde{x}_{k+i} - \bar{x}_{k+i} ||^2_{{Q}} +
        || \tilde{u}_{k+i} - \bar{u}_{k+i} ||^2_{{R}},
\end{split}
\end{align}
where $\bar x$ denotes a reference state, $\bar u$ denotes a reference input, and ${Q_f} \in \mathbb{R}^{n\times n}$, ${Q} \in \mathbb{R}^{n\times n}$, and ${R} \in \mathbb{R}^{m\times m}$ are symmetric positive semi-definite weights. Intuitively, ${Q_f}$ and ${Q}$ place importance on tracking $\bar x_{k:k+N}$, while $R$ penalizes deviations from $\bar u_{k:k+N-1}$.

The performance of NMPC is heavily dependent on the accuracy of the model for $f_d(\cdot)$ in \eqref{eq:eq-constr} to plan optimal trajectories. In practice, NMPC may suffer from unmodeled disturbances and inaccurate estimates of $\theta$. Although this can be remedied through nonlinear parameter estimation, online NMPC can be computationally expensive and may limit the available resources for estimation.

\section{Moving Horizon Parameter Estimation}\label{sec:parameter-estimation}
Similar to NMPC, MHE repeatedly solves an optimization at every time step while respecting hard constraints over a finite time horizon. Given past states, inputs, and a system model, the solution to an MHE problem includes the optimized state, parameter, and disturbance estimates. In this work, we consider a subset of MHE, where the states and process disturbances are already measured or estimated while the parameters are uncertain but bounded. In practice, accurate state and disturbance estimates can be provided by a separate estimation filter.

Assuming that the parameters are not time-varying, we formulate the NMHPE problem as the following NLP:
\begin{subequations}\label{eq:mhe}
\begin{align}
    \min_{ \tilde{\theta} }
        && &V(\tilde \theta)& \label{eq:V} \\
    \text{s.t.}
        && &x_{k-j+1} = f_d({x}_{k-j}, {u}_{k-j}, \tilde{\theta}) + {w}_{k-j}, &\label{nl-constr} \\
        && &\tilde{{\theta}} \in \Theta,& \label{theta-constr}
        % && &{w}_{k-M:k-1} \in \mathcal{W}, \label{w-constr}&
\end{align}
\end{subequations}
where $M \in \mathbb{N}$ is the length of the backward horizon, $j \in \{1,2,\hdots,M\}$ corresponds to a time step within the horizon, and ${w}_{k-M:k-1}$, ${x}_{k-M:k}$, and ${u}_{k-M:k-1}$ are sequences of past disturbances, states, and control inputs, respectively. Constraint \eqref{theta-constr} denotes an empirically determined uncertainty bound on the model parameters.

Since this NMHPE problem attempts to fit parameters to a measured state trajectory, it is natural to treat it as a least squares problem and design $V(\cdot)$ to be quadratic:
\begin{align}
\label{eq:V-quad}
    V(\tilde\theta) :=
        ||\tilde{\theta}-\bar{\theta}||^2_{{P}},
\end{align}
where $\bar\theta$ denotes the previous parameter estimate (i.e., the solutions to \eqref{eq:mhe} at the previous time step), and ${P} \in \mathbb{R}^{p\times p}$ is a symmetric positive semi-definite weighting matrix. In practice, we may also solve for the disturbances in \eqref{eq:mhe} and add a highly-weighted cost term in \eqref{eq:V-quad} to deter deviations from the measured disturbances and guarantee problem feasibility.

\section{Multirotor Dynamics}\label{sec:nl-dynamics}

\subsection{Nonlinear Formulation}
The unperturbed multirotor system considers states in both inertial and body-fixed reference frames. We treat the multirotor UAV as a six degree-of-freedom rigid body, with its state defined as
\begin{align}\label{eq:state}
    x\oft :=& \fourDvecT{p^T\oft}{q^T\oft}{v^T\oft}{\omega^T\oft},
    \shortintertext{where}
    p\oft :=& \threeDvecT{p_x^W\oft}{p_y^W\oft}{p_z^W\oft},\\
    q\oft :=& \fourDvecT{q_w^W\oft}{q_x^W\oft}{q_y^W\oft}{q_z^W\oft},\\
    v\oft :=& \threeDvecT{v_x^B\oft}{v_y^B\oft}{v_z^B\oft},\\
    \omega\oft :=& \threeDvecT{\omega_{x}^B\oft}{\omega_{y}^B\oft}{\omega_{z}^B\oft}.
\end{align}
Here $p\oft$, $q\oft$, $v\oft$, $\omega\oft$, $(\cdot)^W$, and $(\cdot)^B$ denote the vehicle position, attitude, linear velocity, angular velocity, and states relative to the inertial and body-fixed frames, respectively. We elect to parameterize the attitude state as a quaternion to avoid singularities.

\begin{figure}[t!]
    \includegraphics[width=\linewidth]{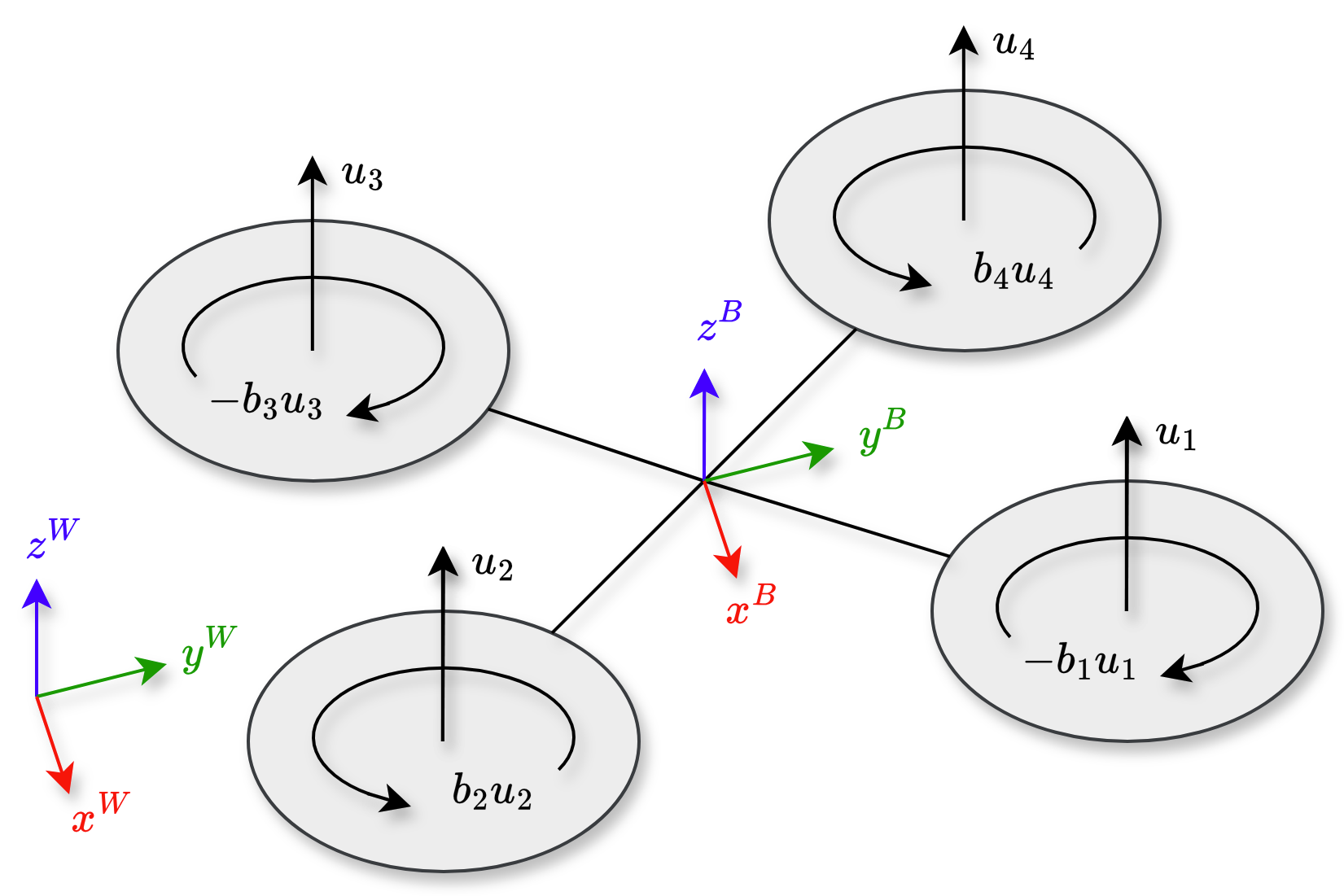}
    \caption{Quadrotor rotor placement, thrusts, and yaw-axis moments relative to the inertial and body-fixed reference frames. Every component of the system's state corresponds to a Cartesian axis in one of the two reference frames.}
\end{figure}

We formulate the unperturbed, continuous-time nonlinear multirotor dynamics \cite{luukkonen2011modelling, RAFFO201029} with the attitude quaternion conventions from \cite{9326337} such that
% \begin{align}\begin{split}\label{eq:nonlinear-dynamics}
%     f(&x\oft, u\oft, \theta) := \\
%     &\begin{bmatrix}
%         \mathcal{Q}(q) v\oft \\
%         \frac{1}{2}G(q)\omega\oft \\
%         \mathcal{Q}^T(q) g + \overm{1}(Ku\oft - A(\theta)v\oft) - \widehat\omega\oft v\oft \\
%         \mathcal{J}^{-1}(\theta)(B(\theta) u\oft - \widehat\omega\oft \mathcal{J}(\theta)\omega\oft)
%     \end{bmatrix},
% \end{split}\end{align}
\begin{align}\begin{split}\label{eq:nonlinear-dynamics}
    f(x, u, \theta) :=
    \begin{bmatrix}
        \mathcal{Q}(q) v \\
        \frac{1}{2}G(q)\omega \\
        \mathcal{Q}^T(q) g + \overm{1}(Ku - A(\theta)v) - \widehat\omega v \\
        \mathcal{J}^{-1}(\theta)(B(\theta) u - \widehat\omega \mathcal{J}(\theta)\omega)
    \end{bmatrix},
\end{split}\end{align}
where $\mathcal{Q}(q)$ is the rotation from the body to the inertial frame, $G(q)$ is the attitude Jacobian, and $\widehat{(\cdot)}$ represents the cross product operator as a skew-symmetric matrix multiplication operation. For the system's physical constants, $\mu$ is the total mass of the system, $g$ is the acceleration due to gravity, $A$ is the diagonal matrix of aerodynamic drag coefficients, and $\mathcal{J}$ is the diagonal inertia matrix:
\begin{align}
    \label{eq:constant_matrices}
    A(\theta) :=& \begin{bmatrix}
        A_{xx} & 0 & 0 \\
        0 & A_{yy} & 0 \\
        0 & 0 & A_{zz}
    \end{bmatrix}, \\
    \mathcal{J}(\theta) :=& \begin{bmatrix}
        I_{xx} & 0 & 0 \\
        0 & I_{yy} & 0 \\
        0 & 0 & I_{zz}
    \end{bmatrix}.
\end{align}
Next, we elect to forgo thrust curve modeling in order to numerically normalize the control inputs. Therefore, the control input $u(t)$ directly consists of the thrust forces produced by each of the drone's rotors. The body frame thrust and torque are defined as $Ku\oft$ and $B(\theta)u\oft$ such that
\begin{align}
    u\oft :=& \fourDvecT{u_{1}\oft}{u_{2}\oft}{\cdots}{u_{m}\oft}, \\
    K :=&
    \begin{bmatrix}
        0_{2 \times m} \\
        1_{1 \times m}
    \end{bmatrix}, \\
    B(\theta) :=&
    \begin{bmatrix}
        d_1 & d_2 & \cdots & d_{m-1} & d_m \\
        -c_1 & -c_2 & \cdots & -c_{m-1} & -c_m \\
        -b_1 & b_2 & \cdots & -b_{m-1} & b_m
    \end{bmatrix},
\end{align}
where $b_i$, $c_i$, and $d_i$ refer to the torque coefficient, body frame horizontal position, and body frame vertical position of rotor $i$, respectively. Finally, we define the parameters of \eqref{eq:nonlinear-dynamics} to be
 \begin{align}\label{eq:parameters}
    \theta :=&
    \begin{bmatrix}
        \mu & \iota^T & a^T & b^T & c^T & d^T
    \end{bmatrix}^T,
    \shortintertext{where}
     \iota :=& \threeDvecT{I_{xx}}{I_{yy}}{I_{zz}}, \\
     a :=& \threeDvecT{A_{xx}}{A_{yy}}{A_{zz}}, \\
     b :=& \fourDvecT{b_1}{b_2}{\cdots}{b_m}, \\
     c :=& \fourDvecT{c_1}{c_2}{\cdots}{c_m}, \\
     d :=& \fourDvecT{d_1}{d_2}{\cdots}{d_m}.
\end{align}

\begin{table*}[htb]
    \centering
    \caption{Crazyflie and Fusion 1 Model Parameters.}
    \label{tab:all-model-parameters}
    \resizebox{0.825\textwidth}{!}{
    \begin{tabular}{| c | c | r | r | c |}
    \hline
    Parameter & Description & Crazyflie & Fusion 1 & Units \\
    \hline
    $\mu$   & Total mass    & $2.70\mathrm{e}-2$ 	& $2.50\mathrm{e}-1$  & kg \\
    $I_{xx}$ 	               & Mass moment of inertia 	 & $1.44\mathrm{e}-5$ 	& $4.27\mathrm{e}-4$ 	& kg m$^2$ \\
    $I_{yy}$ 	               & Mass moment of inertia 	 & $1.40\mathrm{e}-5$ 	& $6.09\mathrm{e}-4$	& kg m$^2$ \\
    $I_{zz}$ 	               & Mass moment of inertia 	 & $2.17\mathrm{e}-5$ 	& $1.50\mathrm{e}-3$  	& kg m$^2$ \\
    $A_{xx}$, $A_{yy}$ 	     & Aerodynamic drag coefficient  & $1.00\mathrm{e}-2$ 	& $2.00\mathrm{e}-2$ 	& kg$/$s \\
    $A_{zz}$ 	               & Aerodynamic drag coefficient 	& $5.00\mathrm{e}-2$ 	& $8.00\mathrm{e}-2$  & kg$/$s \\
    $b_1$, $b_2$, $b_3$, $b_4$ & Yaw-axis torque to thrust ratio  & $2.51\mathrm{e}-2$ 	& $1.11\mathrm{e}-2$    & $-$ \\
    $c_1$, $c_2$ 		       & Body frame horizontal rotor position 	 & $2.83\mathrm{e}-2$ 	& $6.35\mathrm{e}-2$ 	& m \\
    $c_3$, $c_4$ 		       & Body frame horizontal rotor position 	 & $-2.83\mathrm{e}-2$ 	& $-6.35\mathrm{e}-2$	& m \\
    $d_1$, $d_4$ 		       & Body frame vertical rotor position     & $2.83\mathrm{e}-2$ 	& $6.35\mathrm{e}-2$ 	& m \\
    $d_2$, $d_3$ 		       & Body frame vertical rotor position     & $-2.83\mathrm{e}-2$ 	& $-6.35\mathrm{e}-2$	& m \\
    \hline
    \end{tabular}
    }
\end{table*}

\subsection{Affine-in-Parameter Relaxation}\label{sec:aff-dynamics}
The parameters in \eqref{eq:parameters} interact with the system dynamics nonlinearly and are thus difficult to estimate online. Therefore, we are motivated to reformulate the dynamics so that they are affine in the parameters. In \eqref{eq:nonlinear-dynamics}, several of the model parameters are coupled with one another, making linear parameter estimation infeasible. However, by bringing the parameters together into their most compact forms and moving all of the coupled parameters into right-hand side vectors (while still maintaining linearity with respect to the system dynamics), we arrange the second-order terms into the following equations
% \begin{align}
%     \label{eq:prelim-affine1}
%     \dot v &=
%         \mathcal{Q}^T(q) g - \widehat\omega\oft v\oft +
%         Ku\oft/\mu - \mathcal{A}(x\oft)a/\mu, \\
%     \label{eq:prelim-affine2}
%     \dot\omega &=
%         \mathcal{B}(u\oft)\threeDvec{s/I_{xx}}{r/{I_{yy}}}{b/{I_{zz}}} - \mathcal{I}(x\oft)\threeDvec
%             {(I_{zz}-I_{yy})/{I_{xx}}}{(I_{xx}-I_{zz})/{I_{yy}}}{(I_{yy}-I_{xx})/{I_{zz}}},
% \end{align}
\begin{align}
    \label{eq:prelim-affine1}
    \dot v &=
        \mathcal{Q}^T(q) g - \widehat\omega v +
        Ku/\mu - \mathcal{A}(x)a/\mu, \\
    \label{eq:prelim-affine2}
    \dot\omega &=
        \mathcal{B}(u)\threeDvec{d/I_{xx}}{c/{I_{yy}}}{b/{I_{zz}}} - \mathcal{I}(x)\threeDvec
            {(I_{zz}-I_{yy})/{I_{xx}}}{(I_{xx}-I_{zz})/{I_{yy}}}{(I_{yy}-I_{xx})/{I_{zz}}},
\end{align}
where
% \begin{align}
%     \label{eq:block-begin}
%     \mathcal{A}(x\oft) :=& \begin{bmatrix}
%             v_x^B\oft & 0 & 0 \\
%             0 & v_y^B\oft & 0 \\
%             0 & 0 & v_z^B\oft
%         \end{bmatrix}, \\
%     \mathcal{B}(u\oft) :=& \begin{bmatrix}
%             u^T\oft && {0_{1\times2m}} \\
%             0_{1\times m} & -u^T\oft & 0_{1\times m} \\
%             {0_{1\times2m}} && 1^{\mp}_{1\times m} \odot u^T\oft
%         \end{bmatrix}, \\
%     \mathcal{I}(x\oft) :=& \begin{bmatrix}
%             {\omega_y^B\oft \omega_z^B\oft} & 0 & 0 \\
%             0 & {\omega_x^B\oft \omega_z^B\oft} & 0 \\
%             0 & 0 & {\omega_x^B\oft \omega_y^B\oft}
%         \end{bmatrix} \label{eq:block-end}.
% \end{align}
\begin{align}
    \label{eq:block-begin}
    \mathcal{A}(x) :=& \begin{bmatrix}
            v_x^B & 0 & 0 \\
            0 & v_y^B & 0 \\
            0 & 0 & v_z^B
        \end{bmatrix}, \\
    \mathcal{B}(u) :=& \begin{bmatrix}
            u^T && {0_{1\times2m}} \\
            0_{1\times m} & -u^T & 0_{1\times m} \\
            {0_{1\times2m}} && 1^{\mp}_{1\times m} \odot u^T
        \end{bmatrix}, \\
    \mathcal{I}(x) :=& \begin{bmatrix}
            {\omega_y^B \omega_z^B} & 0 & 0 \\
            0 & {\omega_x^B \omega_z^B} & 0 \\
            0 & 0 & {\omega_x^B \omega_y^B}
        \end{bmatrix} \label{eq:block-end}.
\end{align}
The element-wise multiplication operator is denoted by $\odot$, and $(\cdot)^{\mp}$ refers to a vector with entries possessing alternating signs. In this reformulation, by relaxing the parameters and lifting the dynamics, \eqref{eq:parameters} can be defined in new coordinates that \eqref{eq:generic-dynamics} is affine in. We propose a change of variables similar to those of \cite{9290199, doi:10.2514/1.G005376} for all right-hand side coupled parameters in \eqref{eq:prelim-affine1}--\eqref{eq:prelim-affine2} such that the new parameter vector is given by
\begin{align}\label{eq:affine-vec}
    \vartheta :=
    \begin{bmatrix}
        \mathscr{m} &
        \mathscr{a}^T &
        \mathscr{d}^T &
        \mathscr{c}^T &
        \mathscr{b}^T &
        \mathscr{I_{xx}}^T &
        \mathscr{I_{yy}}^T &
        \mathscr{I_{zz}}^T
    \end{bmatrix}^T,
\end{align}
where
\begin{align}
    \mathscr{m} :=& 1/\mu, \\
    \mathscr{a} :=& a/\mu, \\
    \mathscr{d} :=& d/{I_{xx}}, \\
    \mathscr{c} :=& c/{I_{yy}}, \\
    \mathscr{b} :=& b/{I_{zz}}, \\
    \mathscr{I_{xx}} :=& (I_{zz}-I_{yy})/{I_{xx}}, \\
    \mathscr{I_{yy}} :=& (I_{xx}-I_{zz})/{I_{yy}}, \\
    \mathscr{I_{zz}} :=& (I_{yy}-I_{xx})/{I_{zz}}.
\end{align}
By separating the terms in \eqref{eq:prelim-affine1}--\eqref{eq:prelim-affine2} that are multiplied by parameters from those that are not, arranging \eqref{eq:block-begin}-\eqref{eq:block-end} into a linear-in-parameter object, and including \eqref{eq:affine-vec} as the new parameter vector, we obtain the affine-in-parameter dynamics:
% \begin{align}
%     \label{eq:affine-xdot}
%         \dot{x} =& \F + \G\vartheta + w\oft, \\
% \shortintertext{where}
%     \label{eq:F}
%     \F :=& \fourDvec
%         {\mathcal{Q}(q) v\oft}
%         {\frac{1}{2}G(q) \omega\oft}
%         {\mathcal{Q}^T(q) g - \widehat\omega\oft v\oft}
%         {0_{3 \times 1}}, \\
%     \label{eq:G}
%     \G :=& \begin{bmatrix}
%         \multicolumn{3}{c}{0_{7 \times (7+3m)}} \\
%         Ku\oft & -\mathcal{A}(x\oft) & 0_{3 \times (3+3m)} \\
%         0_{3 \times 4} & \mathcal{B}(u\oft) & -\mathcal{I}(x\oft)
%     \end{bmatrix}.
% \end{align}
\begin{align}
    \label{eq:affine-xdot}
        \dot{x} =& \F + \G\vartheta + w, \\
\shortintertext{where}
    \label{eq:F}
    \F :=& \fourDvec
        {\mathcal{Q}(q) v}
        {\frac{1}{2}G(q) \omega}
        {\mathcal{Q}^T(q) g - \widehat\omega v}
        {0_{3 \times 1}}, \\
    \label{eq:G}
    \G :=& \begin{bmatrix}
        \multicolumn{3}{c}{0_{7 \times (7+3m)}} \\
        Ku & -\mathcal{A}(x) & 0_{3 \times (3+3m)} \\
        0_{3 \times 4} & \mathcal{B}(u) & -\mathcal{I}(x)
    \end{bmatrix}.
\end{align}
Despite \eqref{eq:affine-vec} being ``relaxed'', it has the same dimensions as \eqref{eq:parameters}. Although this appears to contradict the intuition that relaxed coordinates should yield more degrees of freedom, \eqref{eq:affine-vec} is indeed less constrained. The lack of increase in dimensionality can be attributed to the fact that several parameters affect the same states and are therefore grouped together. Since some parameters appear multiple times in the decoupled parameter vector, there are new search spaces in regions that were previously nonlinearly constrained.

\section{Implementation and Results}\label{sec:results}

\subsection{Numerical Integration of Nonlinear Dynamics}
To accurately discretize \eqref{eq:nonlinear-dynamics}, we implement fourth order Runge-Kutta (RK4) integration. Furthermore, we normalize the attitude state after every pass of RK4 to prevent the quaternion from drifting away from the four-dimensional unit circle. We utilize this RK4 and quaternion normalization scheme within \eqref{eq:eq-constr}, \eqref{nl-constr}, and the simulated dynamics.

\subsection{Affine-in-Parameter Convexification}
Although RK4 is employed in NMHPE, it is untenable with linear parameter estimation techniques due to its inability to preserve the affine-in-parameter structure of \eqref{eq:affine-xdot} in discrete time. This can be attributed to the nonlinearity of \eqref{eq:affine-xdot} in \eqref{eq:state}, for which the RK4 integration yields multiple components of \eqref{eq:affine-vec} multiplied with one another. To maintain linearity with respect to \eqref{eq:affine-vec}, we instead integrate \eqref{eq:affine-xdot} with the forward Euler method.
Given this forward Euler integration and that \eqref{eq:V-quad} is quadratic with respect to parameters and disturbances, \eqref{eq:mhe} can be convexified into a linear-quadratic form, LQ-MHPE, and solved as a QP. Due to the greater search space of \eqref{eq:affine-vec}, estimating the relaxed parameter vector may lead to dynamically infeasible values. However, in the case of LQ-MHPE, we can transform the uncertainty bounds on the unknown parameters into equivalent constraints on the relaxed parameters, allowing us to tighten the solution search space while still leveraging the numerical benefits of affine-in-parameter dynamics. It is important to note that the original parameters are irrecoverable from the decoupled parameters--the affine-in-parameter formulation exists solely to improve numerical robustness and efficiency of online model adaptation and does not aid in the identification of the true model parameters.

\subsection{Tuning Heuristic}\label{sec:heuristic}
The estimation accuracy of MHPE is highly dependent on the weights of \eqref{eq:V-quad} to numerically scale the estimated parameters. Although the parameters are uncertain, we can make two reasonable assumptions: the nominal multirotor parameters are known, and the nominal parameters exist within known uncertainty bounds. We propose a simple, yet effective heuristic for tuning the MHPE cost weights based on the nominal parameters--set each weight to be a diagonal matrix, and each entry on the diagonal to be the reciprocal of the corresponding parameter's order of magnitude. This heuristic assumes that the uncertainty of the model parameters is proportional to their magnitudes. If this assumption is true, then the heuristic ensures that every parameter is appropriately represented within the numerical optimization while also not requiring any unreasonable prior information.

\begin{figure*}[htb]
    \centering
    \includegraphics[width=\linewidth]{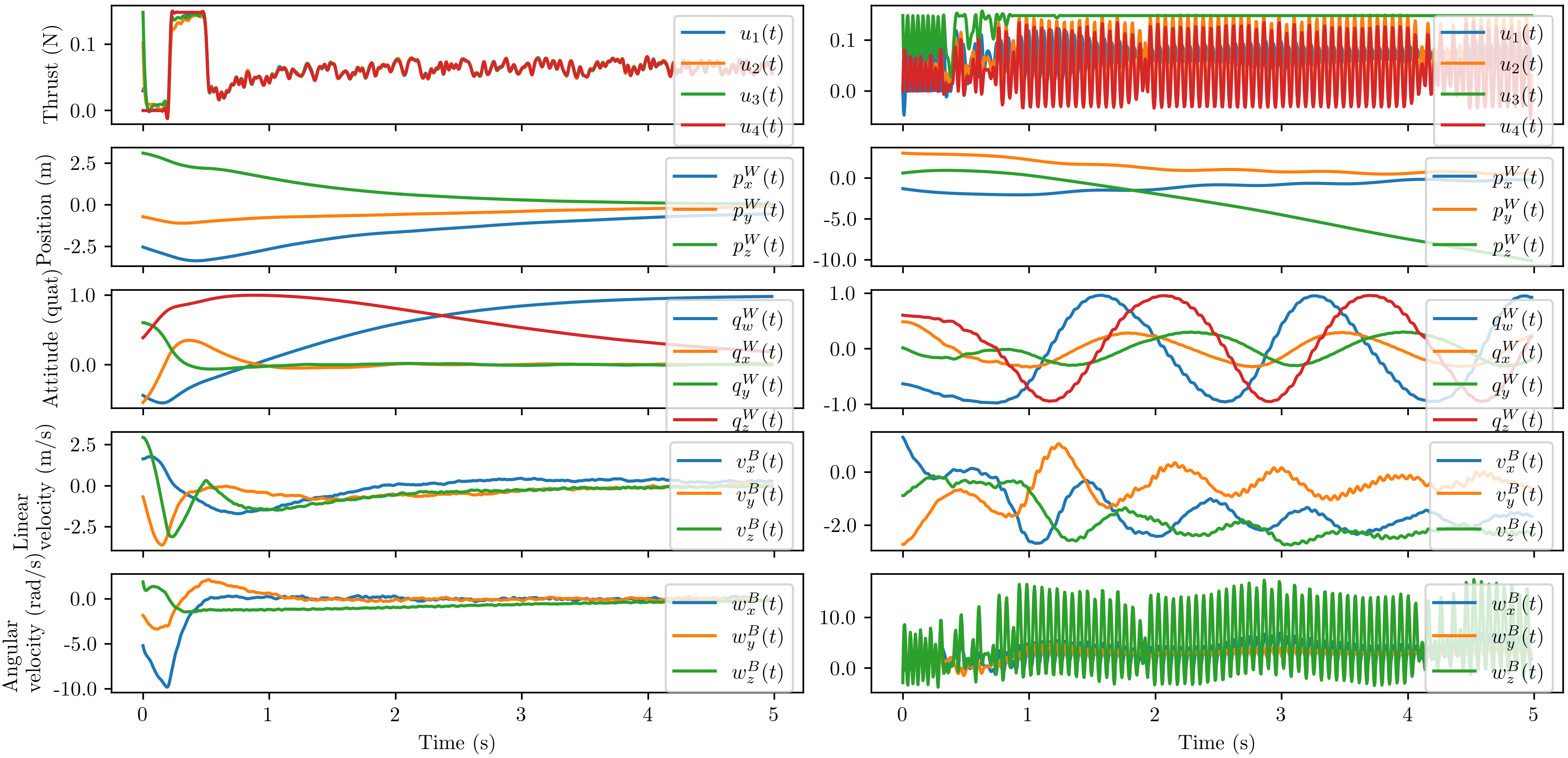}
    \caption{A side-by-side comparison of simulation trials running NMPC. Top to bottom: control inputs, position, attitude, body-fixed linear velocity, and body-fixed angular velocity. Left: an optimal trajectory produced by running NMPC with accurate model parameters. Right: a suboptimal trajectory due to $\pm70\%$ parameter uncertainty within the simulated quadrotor.}
    \label{fig:compare}
\end{figure*}

\subsection{Simulation Environment}
To demonstrate the performance of LQ-MHPE, we benchmark its ability to adapt NMPC against parameter uncertainty. We run Monte Carlo simulations across multiple NMPC control schemes and quadrotor models, where the adaptive NMPC must drive a randomly initialized quadrotor to the origin of the state space. Specifically, we run trials across the nominal Crazyflie \cite{8905295} and Fusion 1 \cite{doi:10.2514/6.2020-1238} quadrotor models. Table \ref{tab:all-model-parameters} contains an overview of each model's parameter values.

Within each trial, the true quadrotor is assigned random model parameters, initialized to a random state, and subject to random disturbances at every time step. All random values are sampled from uniform distributions and subject to bounds on a per-model basis. Table \ref{tab:sim-parameters} contains the bounds in which random values are uniformly sampled.
In every group of trials of the same quadrotor model, trials can be further categorized into having different MHPE schemes, each paired with the same nominal NMPC. Every NMHPE and LQ-MHPE is tuned according to the heuristic described in \ref{sec:heuristic}.
We treat the parameter uncertainty to be uniformly random within prescribed bounds that linearly scale with the given nominal quadrotor's model parameters (e.g., $\theta_{lb} = 0.5\theta_0, \theta_{ub} = 1.5\theta_0$). These bounds are also inputted into MHPE to tighten the parameter estimation search space. For LQ-MHPE, we transform the parameter bounds into the affine-in-parameters space, taking into account the nonlinear effects that the change of variables of \eqref{eq:affine-vec} has on said bounds. See Table \ref{tab:sim-parameters} for the scaling factor uncertainty bounds on the model parameters. We also assume accurate state and disturbance estimates by directly reading the true states and disturbances within the simulation.

The motivation for these trials is to comprehensively explore how LQ-MHPE performs in even the most challenging of control problems. Although these conditions may be more difficult than those typically encountered in real-world drone applications, they serve to capture scenarios in which LQ-MHPE may struggle to adequately perform.

\subsection{Results}
For every test case, that is, every control scheme for every quadrotor model, we simulate 1,000 trials subject to random initial conditions and parameter values. Every trial is simulated for 10 seconds, with 0.02-second time steps. We use \texttt{CasADi}\cite{Andersson2019} to implement NMPC and MHPE, interfacing with \texttt{Ipopt}\cite{Wächter2006} for NMPC and NMHPE and \texttt{OSQP}\cite{osqp} for LQ-MHPE. The convergence tolerances for \texttt{Ipopt} and \texttt{OSQP} are set to 1.0e-8 while their  maximum number of iterations are set to their default values. Although our simulations are run on an AMD Ryzen 7 7840U 5.1 GHz processor that may not be feasible for onboard computing, the relative differences in performance should still carry significance on resource-constrained systems. Additionally, we use NMPC as the baseline because it is the most optimistic, albeit computationally expensive, controller. In practice, a suboptimal linear MPC could be run onboard as a QP to achieve real-time computation.

% Sim Parameter values
\begin{table}[t]
    \centering
    \vspace{0.25cm}
    \caption{Simulation Trial Parameters.}
    \label{tab:sim-parameters}
    \resizebox{0.485\textwidth}{!}{
    \begin{tabular}{| c | r | r |}
    \hline
    Description & Crazyflie & Fusion 1 \\
    \hline
    Parameter lower bound factor & $0.5$ 	& $0.5$ \\
    Parameter upper bound factor & $1.5$ 	& $1.5$ \\
    Process noise lower bound& $-2.5$ 	& $-2.5$ \\
    Process noise upper bound& $2.5$ 	& $2.5$ \\
    Initial position lower bound (m)& $-5.0$ 	& $-10$ \\
    Initial position upper bound (m)& $5.0$ 	& $10$ \\
    Initial linear vel. lower bound (m/s)& $-2.5$ 	& $-5.0$ \\
    Initial linear vel. upper bound (m/s)& $2.5$ 	& $5.0$ \\
    Initial angular vel. lower bound (rad/s)& $-2.5$ 	& $-5.0$ \\
    Initial angular vel. upper bound (rad/s) & $2.5$ 	& $5.0$ \\
    Sampling time (s) & $0.02$ 	& $0.02$ \\
    \hline
    \end{tabular}
    }
\end{table}

% Monte Carlo Stats
\begin{table*}[htb]
    \centering
    \caption{Monte Carlo Simulation Results -- Solve Times and Trajectory Costs.}
    \label{tab:sim-results}
    \resizebox{0.92\textwidth}{!}{
    \begin{tabular}{| c | r r c | r r c |}
    \hline
     & & Crazyflie & & & Fusion 1 & \\
    \hline
    Performance Metric & LQ-MHPE & NMHPE & None & LQ-MHPE & NMHPE & None \\
    \hline
    Best-case solve time (s)& $2.43\mathrm{e}-4$ & $3.34\mathrm{e}-3$ & $-$ & $2.33\mathrm{e}-4$ & $2.96\mathrm{e}-3$ & $-$ \\
    Average solve time (s)& $1.26\mathrm{e}-3$ & $7.02\mathrm{e}-2$ & $-$ & $7.80\mathrm{e}-4$ & $4.32\mathrm{e}-2$ & $-$ \\
    Worst-case solve time (s)& $4.22\mathrm{e}-2$ & $6.51\mathrm{e}-0$ & $-$ & $4.90\mathrm{e}-2$ & $3.92\mathrm{e}-0$ & $-$ \\

    Best-case trajectory cost & $3.59\mathrm{e}+2$ & $3.57\mathrm{e}+2$ & $1.12\mathrm{e}+3$ & $8.01\mathrm{e}+2$ & $8.01\mathrm{e}+2$ & $1.30\mathrm{e}+3$ \\
    Average trajectory cost & $1.28\mathrm{e}+5$ & $2.92\mathrm{e}+5$ & $2.85\mathrm{e}+5$ & $1.21\mathrm{e}+6$ & $1.59\mathrm{e}+6$ & $1.86\mathrm{e}+6$ \\
    Worst-case trajectory cost & $3.52\mathrm{e}+6$ & $4.60\mathrm{e}+6$ & $6.03\mathrm{e}+6$ & $2.26\mathrm{e}+7$ & $2.47\mathrm{e}+7$ & $2.91\mathrm{e}+7$ \\
    \hline
    \end{tabular}
    }
\end{table*}

\begin{figure}[htb]
    \centering
    \includegraphics[width=\linewidth]{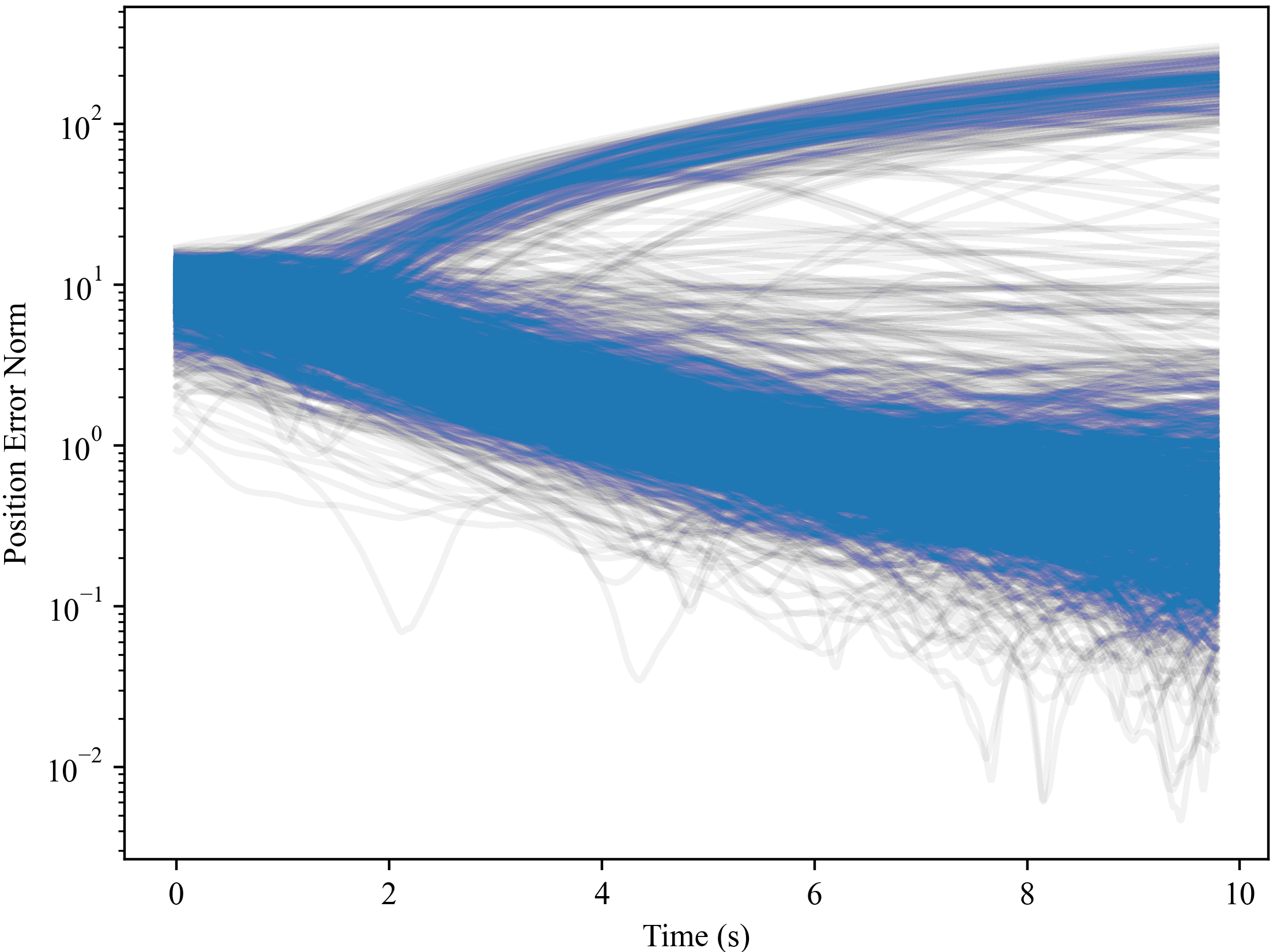}
    \caption{The position error norm trajectories of the NMPC paired with LQ-MHPE on the Fusion 1 quadrotor across 1,000 Monte Carlo simulations. Although some trials have diverging trajectories, there is a strong region of attraction that reduces the position error over time even in the presence of significantly random parameters and initial states.}
    \label{fig:positions}
\end{figure}

\begin{figure}[htb]
    \centering
    \includegraphics[width=\linewidth]{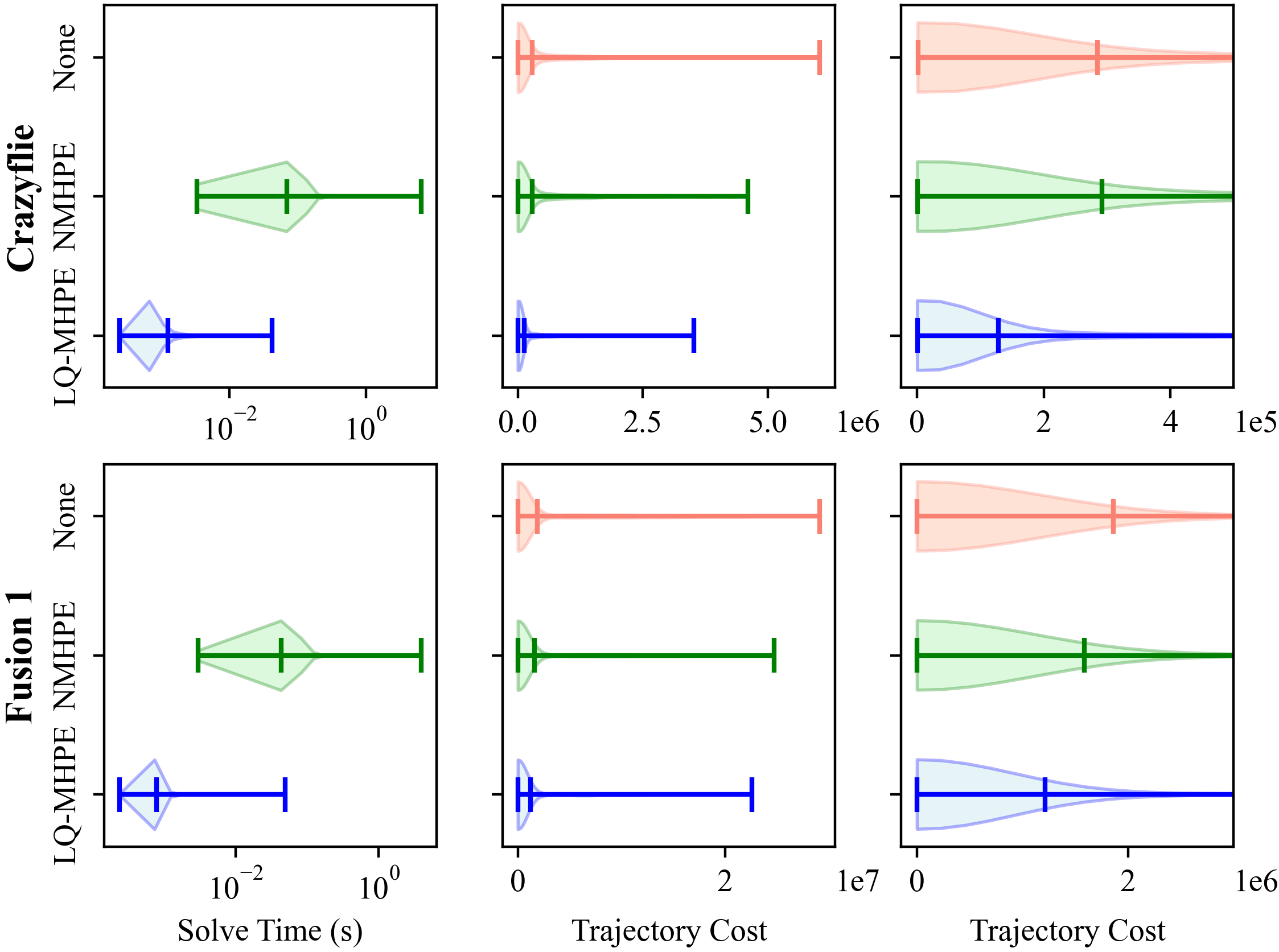}
    \caption{The resulting solve time and trajectory cost distributions for NMPC paired with LQ-MHPE, NMHPE, and no parameter estimator across 1,000 simulation trials for each test case. Plots in the right column zoom in on the trajectory cost distribution means. For all plots, smaller is better.}
    \label{fig:violin}
\end{figure}

Figure~\ref{fig:compare} shows example NMPC simulation trials with and without accurate model parameters, which highlights its inability to handle significant model uncertainty. Figure~\ref{fig:positions} shows the norm of the 3D-position error for all 1,000 trials using NMPC paired with LQ-MHPE converging to a small value for the large majority of simulations, which signifies its robustness to significant model uncertainty. 
Table \ref{tab:sim-results} and Figure \ref{fig:violin} show the resulting MHPE solve times and recorded optimal trajectory costs from the simulation trials. These results show that the affine-in-parameter convexification of the MHPE problem yields significant computational improvements. For the Crazyflie trials, NMHPE had an average solve time of 0.07 seconds. On the other hand, the Crazyflie LQ-MHPE demonstrated an average runtime of around 1 millisecond, indicating a decrease by \textbf{98.2\%}. The Fusion 1 NMHPE and LQ-MHPE performed similarly to their Crazyflie counterparts, with NMHPE and LQ-MHPE yielding 0.0432 and 0.000780 seconds for their mean solve times, respectively. Again, the LQ-MHPE demonstrated an order-of-magnitude decrease in average compute time by \textbf{98.2\%}.

LQ-MHPE also outperforms NMHPE in terms of enhancing NMPC optimality, suggesting that LQ-MHPE is able to produce more accurate parameter estimates. The Crazyflie's mean trajectory cost for the NMPC paired with LQ-MHPE is 1.28e+5, while the NMPC with and without NMHPE demonstrated average costs of 2.92e+5 and 2.85e+5, respectively. In this set of trials, the LQ-MHPE was able to facilitate \textbf{56.2\%} and \textbf{55.1\%} decreases in optimality costs relative to those of the NMPC with and without NMHPE, respectively. Similarly, the LQ-MHPE-based adaptive NMPC decreased trajectory costs by \textbf{23.9\%} and \textbf{34.9\%} compared to those of the NMHPE-based adaptive NMPC and nominal NMPC, respectively. Additionally, LQ-MHPE yielded lower worst-case trajectory costs across both quadrotor models.

One reason why NMHPE computed comparatively poor parameter estimates is because of the highly nonlinear coupling of \eqref{eq:parameters} within \eqref{eq:nonlinear-dynamics}. Although NMHPE should theoretically benefit from the high-fidelity RK4 integration, the nonlinearities associated with \eqref{eq:parameters} increase the difficulty and solve time of the NLP. This may explain why NMHPE-based adaptive NMPC facilitated worse trajectory optimality compared to that of the nominal NMPC in the Crazyflie trials--the poor NLP solutions yielded parameter estimates that degraded NMPC performance and quadrotor stability.

Although LQ-MHPE's forward Euler integration tends to degrade in numerical accuracy compared to that of RK4 over long time steps, the computational speed of LQ-MHPE enables extremely fast (and therefore more accurate) sampling rates. The observed mean compute times for LQ-MHPE are orders of magnitude smaller than the simulated 0.02-second sampling time with minimal tuning of solver parameters, suggesting that these real-time solve rates are also achievable on board when the QP solver is optimized for the processor of choice.

\section{Conclusion and Future Work}\label{sec:conclusion}
The proposed optimal parameter estimator, LQ-MHPE, enables fast, numerically robust adaptive NMPC to facilitate energy-optimal control that adapts to uncertainty in system dynamics. We benchmark its ability to facilitate adaptive NMPC in Monte Carlo simulations, demonstrating a \textbf{98.2\%} decrease in average solve times while reducing trajectory optimality costs by \textbf{23.9-56.2\%}. This framework may enable more real-world use cases of optimization-based parameter estimation in UAV applications where optimal energy consumption and robustness against uncertainties are necessary.

This work can be expanded in multiple directions. The first is onboard hardware implementation--it has been demonstrated that linear MPC can be run in real time on microcontrollers \cite{nguyen2024tinympc}, and the proposed numerical methods should be naturally extendable to LQ-MHPE. The second is to compare LQ-MHPE with simpler linear parameter estimators--such as the Kalman filter--to analyze the trade-offs between computational cost and constraint satisfaction. Finally, there is room to investigate the performance and tuning of the control/parameter estimation scheme in the presence of state estimation (as opposed to idealized state measurements).

\urlstyle{tt}
The source code can be found at \url{https://github.com/EASEL-UCI/parameter-affine-relaxation}.

\footnotesize
\bibliographystyle{IEEEtran}
\bibliography{qrac.bib}

\end{document}